\documentclass[12pt,a4paper]{article}
\usepackage{amsmath}
\usepackage{amssymb}
\usepackage{amsthm}

\usepackage{verbatim}

\usepackage[utf8]{inputenc}%(only for the pdftex engine)
\usepackage{graphicx}
\usepackage{color}
\usepackage{enumerate}

\usepackage{esint}
\usepackage[all]{xy}
\usepackage{framed}
\usepackage{bbm}
\usepackage{subcaption}

\theoremstyle{plain}

\theoremstyle{definition}

\numberwithin{equation}{section}

\newcommand{\bfb}{\mathbf{b}}

\newcommand{\bfa}{\mathbf{a}}
\newcommand{\bfc}{\mathbf{c}}
\newcommand{\bfv}{\mathbf{v}}

\newcommand{\bff}{\mathbf{f}}

\newcommand{\bfe}{\mathbf{e}}

\newcommand{\bfV}{\mathbf{V}}
\newcommand{\bfF}{\mathbf{F}}
\newcommand{\bfw}{\mathbf{w}}

\newcommand{\D}{\mathcal{D}}
\newcommand{\M}{\mathcal{M}}

\newcommand{\dd}{\mathrm{\,d}}

\newcommand{\oO}{\overline{\Omega}}

\newcommand{\mR}{\mathbb{R}}

\DeclareMathOperator{\diag}{diag}

\DeclareMathOperator{\sgn}{sgn}

\DeclareMathOperator{\dist}{dist}

\DeclareMathOperator{\im}{Im}
\DeclareMathOperator{\re}{Re}

\begin{document}

\author{Karl K. Brustad\\ {\small Norwegian University of Science and Technology}}
\title{The one-dimensional nonlocal\\ Dominative $p$-Laplace equation}

\maketitle

\begin{abstract}
The explicit solution to the Dirichlet problem for a class of mean value equations on the real line is derived.
It shed some light on the behavior of solutions to general nonlocal elliptic equations.
\end{abstract}

\section{Introduction}

Imagine you are standing at a point $x$ on a ledge $[-1,1]$.
A number $-1\leq r\leq 1$ is picked at random. If it is positive you get to walk $r$ metres to the right, but if it is negative you have to a take a step $|r|$ metres towards, and perhaps over, the edge at $x=-1$. The process is repeated until you either have reached solid ground at $x>1$, or until you fall into the abyss $x<-1$. What are the chances of surviving?
\begin{figure}[h]
\centering
\includegraphics{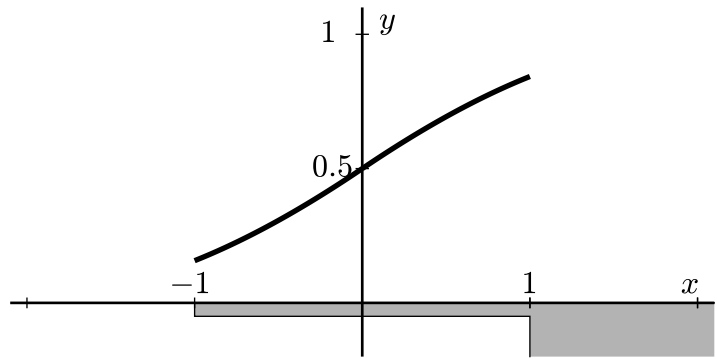}%
\caption{The probability $y$ of surviving when starting from $x\in[-1,1]$.}
\label{fig5}%
\end{figure}

We shall see that the answer is
\[\frac{1}{2}\left(\frac{\cos(1/2)}{1-\sin(1/2)}\sin(x/2) + 1 - \sgn(x)\Big[1 - \cos(x/2)\Big]\right).\]
It is the solution to the mean value equation
\[u(x) = \frac{1}{2}\int_{x-1}^{x+1} u(y)\dd y,\qquad  x\in[-1,1],\]
with boundary conditions
\[u(x) =
\begin{cases}
0,\qquad &-2\leq x < -1,\\
1, & 1< x \leq 2.
\end{cases}\]

The equation above is a special case of the nonlinear problem
\begin{align}
u(x) &= \tfrac{N+2}{N+p}\fint_{B_\epsilon(x)}u(y)\dd y + \tfrac{p-2}{N+p}\sup_{|\xi|=1}\frac{u(x-\epsilon\xi) + u(x+\epsilon\xi)}{2},&& x\in\oO,\label{eq:eq}\\
u(x) &= f(x), && x\in\Gamma^\epsilon,\label{eq:bndcond}
\end{align}
which is investigated in the paper \cite{MR4149517}.
Here, $p\in[2,\infty)$ and $\epsilon>0$ are fixed parameters, $\Omega\subseteq\mR^N$ is open and bounded, and $\Gamma^\epsilon := \{x\notin\oO\;|\; \dist(x,\Omega)\leq \epsilon\}$ is the outer strip of width $\epsilon$. %and $\Omega^\epsilon := \oO\cup \Gamma^\epsilon$.
Furthermore, $f\colon\Gamma^\epsilon\cup\partial\Omega\to\mR$ is a given bounded and integrable function.
\eqref{eq:eq} approximates the Dominative $p$-Laplacian equation
\[0 = \D_p u := \Delta u + (p-2)\lambda_{\max}(D^2u)\]
in the sense that, if we denote the mean value operator by $\M_p^\epsilon$, then
\[\M_p^\epsilon\phi(x) - \phi(x) = \frac{\epsilon^2}{2(N+p)}\D_p\phi(x) + o(\epsilon^2)\]
as $\epsilon\to 0$ for $C^2$ functions $\phi$. The uniformly elliptic operator $\D_p$ was introduced in \cite{MR4085709} in order to explain a superposition principle in the $p$-Laplace equation.

By Lemma 2.2 in \cite{MR4149517} we know that there is a unique solution to \eqref{eq:eq} - \eqref{eq:bndcond}. We are going to derive the solution in the one-dimensional case,
\begin{equation}
\begin{aligned}
u(x) &= \tfrac{3}{p+1}\cdot\tfrac{1}{2\epsilon}\int_{x-\epsilon}^{x+\epsilon} u(y)\dd y + \tfrac{p-2}{p+1}\cdot\tfrac{1}{2}\Big(u(x-\epsilon) + u(x+\epsilon)\Big),\qquad x\in[a,b],\\
u(x) &= f(x),\hspace{195pt} x\in[a-\epsilon,a)\cup(b,b+\epsilon].
\end{aligned}
\label{eq:dirp}
\end{equation}
That is, when $N=1$ and when $\Omega = (a,b)$ is an interval in $\mR$. We shall assume that $\epsilon = (b-a)/n$ for some even number $n=2m$, and that $f$ is continuous.
The main result of \cite{MR4149517} states that the solutions of \eqref{eq:eq} - \eqref{eq:bndcond}, as $\epsilon\to 0$, converge uniformly to the solution of the corresponding local Dirichlet problem. At least for well-behaved domains and boundary values. In our case, this amounts to the simple equation $u''(x) = 0$ and the nonlocal solutions will therefore converge to the affine function with endpoint values $f(a)$ and $f(b)$. Our explicit formulas give some insight to the nature of this convergence. In Section 2 we solve the problem for $p=2$. The case $p>2$, including 
the infinity-equation
\begin{equation}
\begin{aligned}
u(x) &= \frac{1}{2}\Big(u(x-\epsilon) + u(x+\epsilon)\Big), & x&\in[a,b],\\
u(x) &= f(x), & x&\in[a-\epsilon,a)\cup(b,b+\epsilon]
\end{aligned}
\label{eq:infty}
\end{equation}
which is obtained by sending $p\to\infty$, is considered in Section 3.

The stochastic interpretation of the Dirichlet problem \eqref{eq:dirp} is as follows. Suppose you start a random walk from $x_0\in[a,b]$ where each step is chosen from the $\epsilon$-neighbourhood of the previous one according to the rule
\begin{equation}
\begin{cases}
\text{with probability $\frac{3}{p+1}$, the point $x_{k+1}\in [x_k-\epsilon,x_k+\epsilon]$ is picked at random.}\\
\text{with probability $\frac{1}{2}\cdot\frac{p-2}{p+1}$, we set $x_{k+1} = x_k-\epsilon$.}\\
\text{with probability $\frac{1}{2}\cdot\frac{p-2}{p+1}$, we set $x_{k+1} = x_k+\epsilon$.}
\end{cases}
\label{eq:stoch}
\end{equation}
You stop the walk once you have left $[a,b]$ at, say, step $k=\tau$. Then $u(x_0)$ is the \emph{expected value} of the random variable $f(x_\tau)$. In particular, if $f=0$ on $[a-\epsilon,a]$ and $f=1$ on $[b,b+\epsilon]$, then $u(x)$ is the \emph{probability} of exiting at the right when starting the walk from $x$.

Note that the $\sup$ disappears in \eqref{eq:dirp}. Thus, in contrast to the higher dimensional situation, there is no \emph{control} over the stochastic process in one dimension and the equation remains linear for $p>2$.

\section{The nonlocal Laplace equation with uniform distribution}

With $N=1$ and $p=2$ the Dirichlet problem reads
\begin{equation}
\begin{cases}
u(x) = \displaystyle\frac{1}{2\epsilon}\int_{x-\epsilon}^{x+\epsilon} u(y)\dd y,\qquad & \text{for $x\in[a,b]$},\\
u(x) = f(x), & \text{for $x\in[a-\epsilon,a)\cup(b,b+\epsilon]$}.
\end{cases}
\label{eq:s}
\end{equation}
If $f$ is constant on each of the two boundary parts, we shall show that $u$ is a \emph{piecewise trigonometric function}. Specifically, on each of the $n$ intervals $[a + (k-1)\epsilon, a + k\epsilon]$, $k=1,\dots,n$, of length $\epsilon$, $u$ is on the form
\begin{equation}
u(x) = a_k + \sum_{j=1}^m b_{k,j}\sin\left(\frac{\lambda_j}{\epsilon}x\right) + c_{k,j}\cos\left(\frac{\lambda_j}{\epsilon}x\right)
\label{eq:uform}
\end{equation}
for computable coefficients $a_k, b_{k,j}, c_{k,j}$, and where
\[\lambda_j = \cos\left(\frac{j\pi}{n+1}\right).\]
When $f$ is not constant, the above formula is supplemented with some additional terms involving trigonometric convolutions of the data. These terms are written out in \eqref{eq:addterms} for the case $[a,b]=[-1,1]$.
\begin{figure}[h]
\centering
\includegraphics{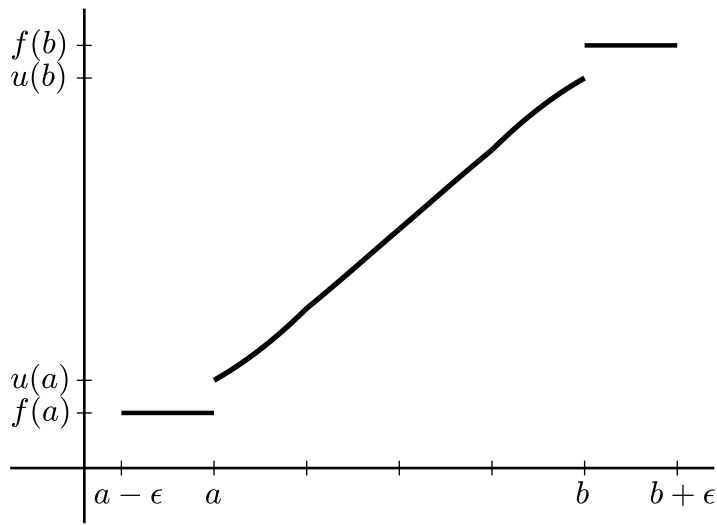}%
\caption{The solution of \eqref{eq:s} with constant boundary values and $n=4$.}%
\label{fig1}%
\end{figure}
\begin{figure}[h]
\centering
\includegraphics{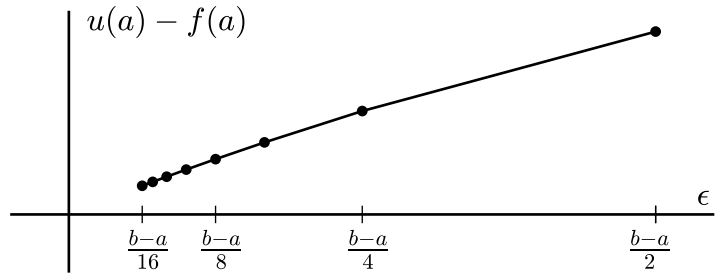}%
\caption{Convergence of the boundary values.}%
\label{fig3}%
\end{figure}

Although the solution has oscillations for all $\epsilon = (b-a)/n$, our plots show that the graph is almost indistinguishable from a straight line on the inner part of the interval already from $n\geq 4$. However, near the endpoints the graph is unmistakably curved and the numerics indicate that the convergence $u(a)\to f(a)$ and $u(b)\to f(b)$ is no better than linear in $\epsilon$. See Figure \ref{fig3}. On the other hand, the comparison principle (Lemma 2.5 \cite{MR4149517}) ensures that the convergence \emph{is} linear, since the solution has to lie between the two lines passing through the points $(a-\epsilon,f(a)), (b,f(b))$ and
$(a,f(a)), (b+\epsilon,f(b))$, respectively. Unfortunately, this argument is purely one-dimensional and does not work for $N\geq 2$.

\subsection{The functional differential equation}

One may easily show that the solution $u$ of \eqref{eq:s} is $C/(2\epsilon)$-Lipschitz on $[a,b]$, where $C = \max f - \min f$. It follows that $u$ is differentiable in $(a,b)\setminus\{a+\epsilon,b-\epsilon\}$ with
\begin{equation}
u'(x) = \frac{1}{2\epsilon}\Big( u(x+\epsilon) - u(x-\epsilon)\Big).
\label{eq:der}
\end{equation}
$u$ is generally not differentiable at $a+\epsilon$ or $b-\epsilon$ since it is no reason $u$ should be continuous at $a$ or $b$. A close inspection of Figure \ref{fig1} reveals kinks in the graph of $u$ at these points.
The formula \eqref{eq:der} provides, however, a weak derivative for $u$ on $[a,b]$.

\subsection{Constant boundary values}

Let $u$ be the solution of
\begin{equation}
\begin{cases}
u(x) = \displaystyle\frac{1}{2\epsilon}\int_{x-\epsilon}^{x+\epsilon} u(y)\dd y,\qquad & \text{for $x\in[a,b]$},\\
u(x) = c_l, & \text{for $x\in[a-\epsilon,a)$},\\
u(x) = c_r, & \text{for $x\in(b,b+\epsilon]$}.
\end{cases}
\label{eq:s2}
\end{equation}
where $a<b$, $\epsilon = (b-a)/n$, for some even number $n=2m$, and where $c_l\neq c_r$ are two constants. If $c_l = c_r$ the solution is identically equal to this common constant.

To exploit the symmetry in the problem, we assume that $u$ is scaled, shifted, and translated so that $[a,b] = [-1,1]$, and $c_l = -1$ and $c_r = 1$. The solution will then be odd. This is possible because the equation is linear and translation invariant. Also, constants are solutions. 

Divide the domain into $n$ parts of length $\epsilon$.
For $k=1,\dots,n$ define the functions $v_k\colon[0,1]\to\mR$ as
\begin{equation}
\begin{aligned}
v_k(t) &:= \frac{2}{c_r-c_l}\left(u\big(\epsilon t +a + (k-1)\epsilon\big) - \frac{c_l+c_r}{2}\right)\\
       &= u\big(\epsilon t - 1 + (k-1)\epsilon\big).
\end{aligned}
\label{eq:trans}
\end{equation}
Now, each $v_k$ is differentiable in $(0,1)$. For $k=1$ we have by \eqref{eq:der}
\begin{align*}
v_1'(t)
	&= \epsilon u'(\epsilon t - 1)\\
	&= \epsilon\frac{1}{2\epsilon}\Big( u(\epsilon t - 1 +\epsilon) - u(\epsilon t - 1 -\epsilon)\Big)\\
	&= \frac{1}{2}\Big( u(\epsilon t - 1 +\epsilon) - c_l\Big)\\
	&= \frac{1}{2}\Big( v_2(t) + 1\Big).
\end{align*}
Similarly, $v_n'(t) = \frac{1}{2}\Big( 1 - v_{n-1}(t)\Big)$, and
for $k = 2,\dots,n-1$ we simply have
\[v_k'(t) = \frac{1}{2}\Big( v_{k+1}(t) - v_{k-1}(t)\Big).\]
This defines a non-homogeneous linear system of ODEs,
\[\begin{bmatrix}
	v'_1\\ v_2'\\ v_3'\\ \vdots\\ v_{n-1}'\\ v_n'
\end{bmatrix} = \frac{1}{2}
\begin{bmatrix}
	0 & 1 & 0 & 0 & \cdots & 0\\
	-1 & 0 & 1 & 0 & \cdots & 0\\
	0 & -1 & 0 & 1 & \cdots & 0\\
	&& \vdots &&&\\
	0 & \cdots & 0 & -1 & 0 & 1\\
	0 & \cdots & 0 & 0 & -1 & 0
\end{bmatrix}
\begin{bmatrix}
	v_1\\ v_2\\ v_3\\ \vdots\\ v_{n-1}\\ v_n
\end{bmatrix} + \frac{1}{2}\begin{bmatrix}
	1\\ 0\\ 0\\ \vdots\\ 0\\ 1
\end{bmatrix},\]
or
\begin{equation}
\bfv'(t) = A\bfv(t) + \bfc
\label{eq:ode}
\end{equation}
in vector notation.
The general solution is $\bfv(t) = e^{tA}(\bfv_0 + A^{-1}\bfc) - A^{-1}\bfc$.
where
\[e^A := \sum_{k=0}^\infty\frac{1}{k!}A^k\]
is the matrix exponential.

Defining
\[\tilde{\bfe} := \sum_{k=1}^n (-1)^k\bfe_k,\]
we note that
\[A\tilde{\bfe} = \frac{1}{2}
\begin{bmatrix}
	0 & 1 & 0 & 0 & \cdots & 0\\
	-1 & 0 & 1 & 0 & \cdots & 0\\
	0 & -1 & 0 & 1 & \cdots & 0\\
	& \vdots &&& \cdots & 1\\
	0 & \cdots && 0 & -1 & 0
\end{bmatrix}\begin{bmatrix}
	-1\\ 1\\ -1\\ \vdots\\ 1
\end{bmatrix} = \frac{1}{2}\begin{bmatrix}
	1\\ 0\\ \vdots\\ 0 \\ 1
\end{bmatrix} = \bfc\]
and thus
\begin{equation}
\bfv(t) = e^{tA}(\bfv_0 + \tilde{\bfe}) - \tilde{\bfe}.
\label{eq:1}
\end{equation}

The initial value $\bfv_0 = \bfv(0)$ contains the values of $u$ at the nodes $-1 + (k-1)\epsilon$ and is of course unknown. By definition, we have the $n-1$ identities $v_k(0) = v_{k-1}(1)$, $k = 2,\dots,n$. An $n$'th equation can be obtained by using the fact that $u$ is an odd function. For example, $u(-1) = - u(1)$ which corresponds to 
$v_1(0) = -v_n(1)$. Thus,
\begin{equation}
\bfv(0) = B\bfv(1)
\label{eq:2}
\end{equation}
where $B$ is the orthogonal matrix
\[B := \begin{bmatrix}
	0 & 0 & 0 & \cdots & 0 & -1\\
	1 & 0 & 0 & \cdots & 0 & 0\\
	0 & 1 & 0 & \cdots & 0 & 0\\
	&& \ddots &&&\\
	0 & \cdots && 0 & 1 & 0
\end{bmatrix} = \begin{bmatrix}
	\bf0^\top & -1\\ I_{n-1} & \bf0
\end{bmatrix}.\]
Inserting \eqref{eq:1} into \eqref{eq:2} with $t=1$, the continuity of $u$ in $[-1,1]$ gives a linear equation for $\bfv_0$,
\[(I-Be^A)\bfv_0 = B(e^A-I)\tilde{\bfe},\]
and the unique solution of \eqref{eq:ode}, \eqref{eq:2} is, after some simplifications,
\[\bfv(t) = e^{tA}(I-Be^A)^{-1}(I-B)\tilde{\bfe} - \tilde{\bfe},\qquad 0\leq t\leq 1.\]
The invertability of $I - Be^A$ is discussed later.
Observe that $\bfv(t)$ does not depend on anything but $n$.

The solution $u$ of the original problem \eqref{eq:s2} can now be assembled by inverting the relations in \eqref{eq:trans}:

\begin{align*}
u(x) &=
\begin{cases}
c_l,\qquad &x\in[a-\epsilon, a),\\
\frac{c_r-c_l}{2}v_k\left(\frac{x-a}{\epsilon} + 1 - k\right) + \frac{c_r+c_l}{2}, &x\in[a + (k-1)\epsilon, a + k\epsilon],\quad k = 1,\dots,n,\\
c_r,\qquad &x\in(b,b+\epsilon].
\end{cases}\notag \\
			&=
\begin{cases}
-1,\qquad &x\in[-1-\epsilon, -1),\\
v_k\left(\frac{x+1}{\epsilon} + 1 - k\right), &x\in[-1 + (k-1)\epsilon, -1 + k\epsilon],\quad k = 1,\dots,n,\\
1,\qquad &x\in(1,1+\epsilon].
\end{cases}
\end{align*}

\subsection{Analysis}

The coefficient matrix of the linear system is the skew-symmetric and tridiagonal $n\times n$ matrix
\[A := \frac{1}{2}
\begin{bmatrix}
	0 & 1 & 0 & 0 & \cdots & 0\\
	-1 & 0 & 1 & 0 & \cdots & 0\\
	0 & -1 & 0 & 1 & \cdots & 0\\
	& \vdots &&& \cdots & 1\\
	0 & \cdots && 0 & -1 & 0
\end{bmatrix}\]
where $n=2m$ is even.

In order to diagonalize $A$, we define the vectors $\xi_j\in\mathbb{C}^n$ with components
\[\bfe_k^\top \xi_j = c_n i^{k+2j}\sin\left(kj\frac{\pi}{n+1}\right),\qquad c_n := \sqrt{\frac{2}{n+1}},\; i :=\sqrt{-1}.\]
They have unit length since
\begin{align*}
\sum_{k=1}^n\sin^2\left(kj\frac{\pi}{n+1}\right)
	&= \frac{n+1}{2}
\end{align*}
by Lagrange's identity. Using the rule
\[\sin\theta + \sin\phi = 2\cos\left(\frac{\theta-\phi}{2}\right)\sin\left(\frac{\theta+\phi}{2}\right),\]
shows that
\begin{align*}
\bfe_k^\top A\xi_j
	&= \frac{c_n}{2}\left(i^{k+2j+1}\sin\left((k+1)j\frac{\pi}{n+1}\right) - i^{k-1+2j}\sin\left((k-1)j\frac{\pi}{n+1}\right) \right)\\
	&= \frac{i^{k+1+2j}c_n}{2}\left(\sin\left((k+1)j\frac{\pi}{n+1}\right) + \sin\left((k-1)j\frac{\pi}{n+1}\right) \right)\\
	&= i^{k+1+2j}c_n\cos\left(j\frac{\pi}{n+1}\right)\sin\left(kj\frac{\pi}{n+1}\right)\\
	&= i\cos\left(j\frac{\pi}{n+1}\right)\bfe_k^\top\xi_j.
\end{align*}
That is, $\xi_j$ is an eigenvector to $A$ with eigenvalue $i\cos\left(j\frac{\pi}{n+1}\right)$. Since the eigenvalues are distinct, this also implies that the eigenvectors are orthogonal. Moreover, for $j=1,\dots,m$ one can show that $\xi_{n+1-j} = \overline{\xi}_j$,
and $A$ is thus diagonalized by the unitary matrix
\[U := [\xi_1, \overline{\xi}_1,\dots, \xi_m, \overline{\xi}_m] \in \mathbb{C}^{n\times n},\]
producing
\[\Lambda := \overline{U}^\top A U = \diag(i\lambda_1,-i\lambda_1,\dots,i\lambda_m,-i\lambda_m)\]
where
\[0 < \lambda_j := \cos\left(\frac{j}{n+1}\pi\right) < 1,\qquad j=1,\dots,m.\]

The real and imaginary parts of the eigenvectors $\xi_j = \bfa_j + i\bfb_j$ are
\begin{align*}
\bfa_j &= c_n\sum_{k=1}^m (-1)^{k+j}\sin\left(2kj\frac{\pi}{n+1}\right)\bfe_{2k},\\
\bfb_j &= c_n\sum_{k=1}^m (-1)^{k+j+1}\sin\left((2k-1)j\frac{\pi}{n+1}\right)\bfe_{2k-1}.
\end{align*}
Here,
\[\bfa_j^\top \bfb_k = 0\qquad\text{and}\qquad \bfa_j^\top \bfa_k = \bfb_j^\top \bfb_k = \frac{\delta_{j,k}}{2}\]
since $\xi_j^\top\xi_k = 0$ and $\overline{\xi_j}^\top\xi_k = \delta_{j,k}$.
If we define the real skew-symmetric matrices
\[A_j := -2\im \xi_j\overline{\xi}_j^\top = 2(\bfa_j\bfb_j^\top-\bfb_j\bfa_j^\top),\qquad j=1,\dots,m,\]
then
\begin{align*}
A &= U\Lambda\overline{U}^\top\\
  &= \sum_{j=1}^m i\lambda_j\left(\xi_j\overline{\xi}_j^\top - \overline{\xi}_j\xi_j^\top\right)\\
	&= \sum_{j=1}^m \lambda_j A_j.
\end{align*}
Furthermore,
\[A_j^2 = -2(\bfa_j\bfa_j^\top + \bfb_j\bfb_j^\top) = -2\re \xi_j\overline{\xi}_j^\top\]
is the negative of a two-rank symmetric projection and
\begin{align*}
e^{tA} &= Ue^{t\Lambda}\overline{U}^\top\\
  &= \sum_{j=1}^m e^{i\lambda_j t}\xi_j\overline{\xi}_j^\top + e^{-i\lambda_j t}\overline{\xi}_j\xi_j^\top\\
	&= 2\re\sum_{j=1}^m e^{i\lambda_j t}\xi_j\overline{\xi}_j^\top\\
	&= \sum_{j=1}^m \sin(\lambda_j t)A_j - \cos(\lambda_j t)A_j^2.
\end{align*}

The solution on the interval $[-1 + (k-1)\epsilon, -1 + k\epsilon]$ is therefore
\begin{align*}
u(x) &= v_k\left(\frac{x+1}{\epsilon} + 1-k\right)\\
     &= \bfe_k^\top\bfv\left(\frac{x+1}{\epsilon} + 1-k\right)\\
		 &= \bfe_k^\top e^{\left(\frac{x+1}{\epsilon} + 1-k\right)A}(\bfv_0 + \tilde{\bfe}) - \bfe_k^\top\tilde{\bfe}\\
		 &= (-1)^{k+1} + \bfe_k^\top e^{(m+1-k)A}e^{mxA}(\bfv_0 + \tilde{\bfe}),\qquad \frac{1}{\epsilon} = m,\\
		 &= (-1)^{k+1} + \bfe_k^\top e^{(m+1-k)A}\left(\sum_{j=1}^m \sin(\lambda_jmx)A_j - \cos(\lambda_jmx)A_j^2\right)(\bfv_0 + \tilde{\bfe}),
\end{align*}
which is on the form \eqref{eq:uform}.

\subsection{The case $n=2$}

When $n=2$,
\[A = \frac{1}{2}\begin{bmatrix}
	0 & 1\\ -1 & 0
\end{bmatrix}\]
with eigenvalues $\pm i/2$. We have $\lambda_1 = 1/2 = \cos(1\cdot\pi/(n+1))$ and 
\[A_1 = \begin{bmatrix}
	0 & 1\\ -1 & 0
\end{bmatrix}\]
so that $A = \lambda_1 A_1$. Next, $A_1^2 = -I$ and it follows that
\[e^{tA} = \sin(\lambda_1 t)A_1 + \cos(\lambda_1 t)I = \begin{bmatrix}
	\cos(t/2) & \sin(t/2)\\ -\sin(t/2) & \cos(t/2)
\end{bmatrix}.\]
Also, $B = -A_1 = A_1^\top$ and
\[I - Be^A = (1 - \sin(1/2))I + \cos(1/2)A_1 = \begin{bmatrix}
	1 - \sin(1/2)& \cos(1/2)\\ -\cos(1/2) & 1 - \sin(1/2)
\end{bmatrix}\]
with inverse
\[(I - Be^A)^{-1} = \frac{1}{2(1-s)}\begin{bmatrix}
	1 - s& -c\\ c & 1 - s
\end{bmatrix} = \frac{(I - Be^A)^\top}{2(1-s)} = \frac{I - e^{-A}A_1}{2(1-s)}.\]
This makes
\begin{align*}
\bfv_0 &= (I - Be^A)^{-1}(Be^A-B)\tilde{\bfe}\\
       &= \frac{1}{2(1-s)}(I - A_1e^{-A})(I-e^A)A_1\tilde{\bfe}\\
			 &= \frac{1}{2(1-s)}\left(I - sA_1 - cI - A_1(-sA_1 + cI) + A_1\right)A_1\tilde{\bfe}\\
			 &= \frac{1-\sin(1/2)-\cos(1/2)}{1-\sin(1/2)}\begin{bmatrix}
					1\\ 0
				\end{bmatrix},
\end{align*}
and
\begin{align*}
\bfv(t) &= e^{tA}(\bfv_0 + \tilde{\bfe}) - \tilde{\bfe}\\
        &= \begin{bmatrix}
	\cos(t/2) & \sin(t/2)\\ -\sin(t/2) & \cos(t/2)
\end{bmatrix}\begin{bmatrix}
					-\frac{\cos(1/2)}{1-\sin(1/2)}\\ 1
				\end{bmatrix} - \begin{bmatrix}
					-1\\ 1
				\end{bmatrix}.
\end{align*}
For $k = 1,2$, the solution in $[-2 + k, -1 + k]$  is therefore $u(x) = v_k\left(x+2 - k\right)$. That is,
\begin{align*}\label{eq:n2}
u(x)	&=
\begin{cases}
-1, &x\in[-2, -1),\\
-C\cos(x/2+1/2) + \sin(x/2+1/2) + 1, \qquad &x\in[-1,0],\\
C\sin(x/2) + \cos(x/2) - 1, &x\in[0,1]\\
1,\qquad &x\in(1,2],
\end{cases}\\
      &=
\begin{cases}
-1, &x\in[-2, -1),\\
C\sin(x/2) - \sgn(x)\Big[1-\cos(x/2)\Big], \hspace{52pt} &x\in[-1,1],\\
1,\qquad &x\in(1,2],
\end{cases}
\end{align*}
where
\[C := \frac{\cos(1/2)}{1-\sin(1/2)}.\]
It is the graph of $(u+1)/2$ that is shown in Figure \ref{fig5}.

\subsection{Non-constant $f$}

Let $u$ be the solution of
\[
\begin{cases}
u(x) = \displaystyle\frac{1}{2\epsilon}\int_{x-\epsilon}^{x+\epsilon} u(y)\dd y,\qquad & \text{for $x\in[-1,1]$},\\
u(x) = f(x), & \text{for $x\in[-1-\epsilon,-1)\cup (1,1+\epsilon]$},
\end{cases}
\]
where $\epsilon = 2/n$. Translate the data to the interval $[0,1]$ by writing
\[f_l(t) := f(\epsilon t - 1 -\epsilon),\qquad f_r(t) := f(\epsilon t+1),\]
and, as usual,
define the functions $v_k\colon[0,1]\to\mR$ as
\[v_k(t) = u(x_k + \epsilon t),\qquad x_k :=  - 1 + (k-1)\epsilon.\]
As before,
\[v_k'(t) = \frac{1}{2}\Big( v_{k+1}(t) - v_{k-1}(t)\Big)\]
for $k = 2,\dots,n-1$, and now
\begin{align*}
v_1'(t)
	&= \frac{1}{2}\Big( v_2(t) - f_l(t)\Big),\\
v_n'(t)
	&= \frac{1}{2}\Big( f_r(t) - v_{n-1}(t)\Big).
\end{align*}
The system then reads
\[\begin{bmatrix}
	v'_1\\ v_2'\\ v_3'\\ \vdots\\ v_{n-1}'\\ v_n'
\end{bmatrix} = \frac{1}{2}
\begin{bmatrix}
	0 & 1 & 0 & 0 & \cdots & 0\\
	-1 & 0 & 1 & 0 & \cdots & 0\\
	0 & -1 & 0 & 1 & \cdots & 0\\
	&& \vdots &&&\\
	0 & \cdots & 0 & -1 & 0 & 1\\
	0 & \cdots & 0 & 0 & -1 & 0
\end{bmatrix}
\begin{bmatrix}
	v_1\\ v_2\\ v_3\\ \vdots\\ v_{n-1}\\ v_n
\end{bmatrix} + \frac{1}{2}\begin{bmatrix}
	-f_l\\ 0\\ 0\\ \vdots\\ 0\\ f_r
\end{bmatrix},\]
or
\[\bfv'(t) = A\bfv(t) + \bff(t)\]
in obvious notation. The general solution is
\begin{equation}
\bfv(t) = e^{tA}\left(\bfv_0 + \int_0^t e^{-sA}\bff(s)\dd s\right), \qquad 0\leq t\leq 1,
\label{eq:gensol}
\end{equation}
and the solution $u$ on the $k$'th interval $[x_k, x_{k+1}]$ is
\begin{align*}
u(x) &= \bfe_k^\top\bfv\left(\frac{x+1}{\epsilon} + 1 - k\right)\\
     &= \bfe_k^\top e^{\left(\frac{x+1}{\epsilon} + 1 - k\right)A}\left(\bfv_0 + \int_0^{\frac{x+1}{\epsilon} + 1 - k} e^{-sA}\bff(s)\dd s\right)\\
		 &= \bfe_k^\top e^{\left(\frac{x+1}{\epsilon} + 1 - k\right)A}\bfv_0 + \frac{1}{\epsilon}\int_{-1+(k-1)\epsilon}^{x} \bfe_k^\top e^{\frac{x-y}{\epsilon}A}\bff\left(\frac{y+1}{\epsilon}+1-k\right)\dd y.
\end{align*}
The first term is on the form \eqref{eq:uform} while the integral may be expanded as
\begin{align}\label{eq:addterms}
\frac{1}{2\epsilon}\sum_{j=1}^m & \int_{-1+(k-1)\epsilon}^{x}  \bigg[\sin\left(\lambda_j(x-y)/\epsilon\right) \Big( a_{k,j}f\big(y+2+(1-k)\epsilon\big) - b_{k,j}f\big(y-k\epsilon\big) \Big)\notag\\
	& {} - \cos\left(\lambda_j(x-y)/\epsilon\right) \Big( c_{k,j}f\big(y+2+(1-k)\epsilon\big) - d_{k,j}f\big(y-k\epsilon\big) \Big) \bigg]\dd y
\end{align}

In order to find an expression for the initial condition $\bfv_0 = \bfv(0)$ we still have the $n-1$ identities
\[v_k(0) = v_{k-1}(1),\qquad k = 2,\dots,n,\]
but now we cannot assume $u$ to be odd and $v_1(0) = - v_n(1)$ is therefore not valid. Thus we are one equation short of determining $\bfv_0$.

Set
\[F_0 := \fint_{\Gamma^\epsilon}f\dd x\]
to be the average of the data and let $U_0 := \int_{-1}^1 u\dd x$ be the total integral of the solution. Then
\begin{align*}
\sum_{k=1}^n v_k(0) + v_n(1)
	&= \sum_{k=0}^n u(-1+ k\epsilon)\\
	&= \frac{1}{2\epsilon}\sum_{k=0}^n \int_{-1 + (k-1)\epsilon}^{-1+(k+1)\epsilon}u(y)\dd y\\
	&= \frac{1}{2\epsilon}\left(\int_{-1-\epsilon}^{-1} f(y)\dd y + 2\int_{-1}^1 u(y)\dd y + \int_1^{1+\epsilon} f(y)\dd y\right)\\
	&= F_0 + \frac{1}{\epsilon}U_0,
\end{align*}
while
\begin{align*}
\sum_{k=1}^m v_{2k-1}(0) + v_n(1)
	&= \sum_{k=0}^m u(-1 + 2k\epsilon)\\
	&= \frac{1}{2\epsilon}\sum_{k=0}^m \int_{-1+(2k-1)\epsilon}^{-1+(2k+1)\epsilon}u(y)\dd y\\
	&= F_0 + \frac{1}{2\epsilon}U_0.
\end{align*}
It follows that
\begin{align*}
F_0 &= 2\left(\sum_{k=1}^m v_{2k-1}(0) + v_n(1)\right) - \left(\sum_{k=1}^n v_k(0) + v_n(1)\right)\\
		&= \sum_{k=1}^n (-1)^{k-1}v_k(0) + v_n(1)\\
		&= v_1(0) + \sum_{k=1}^n (-1)^k v_k(1).
\end{align*}
We now have $n$ equations for $\bfv_0$, namely
\[\begin{bmatrix}
	v_1(0)\\ v_2(0)\\ v_3(0)\\ \vdots\\ v_n(0)
\end{bmatrix} = 
\begin{bmatrix}
	1 & -1 & 1 & -1 & \cdots & -1\\
	1 & 0 & 0 & 0 & \cdots & 0\\
	0 & 1 & 0 & 0 & \cdots & 0\\
	&& \vdots &&&\\
	0 & \cdots && 0 & 1 & 0
\end{bmatrix}
\begin{bmatrix}
	v_1(1)\\ v_2(1)\\ v_3(1)\\ \vdots\\ v_n(1)
\end{bmatrix} + \begin{bmatrix}
	F_0\\ 0\\ 0\\ \vdots\\ 0
\end{bmatrix},\]
or
\[\bfv_0 = \tilde{B}\bfv(1) + F_0\bfe_1.\]
By the general solution \eqref{eq:gensol} we may write
\begin{align*}
\bfv_0 &= \tilde{B}\bfv(1) + F_0\bfe_1\\
	&= \tilde{B}e^{A}\left(\bfv_0 + \int_0^1 e^{-sA}\bff(s)\dd s\right) + F_0\bfe_1,
\end{align*}
which means that
\[\bfv_0 = (I-\tilde{B}e^{A})^{-1}\left(\tilde{B}e^{A}\int_0^1 e^{-sA}\bff(s)\dd s + F_0\bfe_1\right).\]

\begin{figure}[h]
\centering
\includegraphics{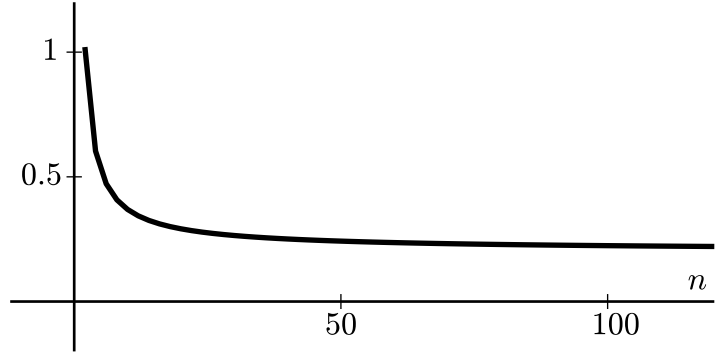}%
\caption{The graph of $\det^{1/n}(I-Be^A)$ as a function of $n$. The corresponding picture with $\tilde{B}$ instead of $B$ is rather similar.}%
\label{fig4}%
\end{figure}

It seems hard to prove that $I-Be^A$ and $I-\tilde{B}e^A$ are invertible or, equivalently, that $Be^A$ and $\tilde{B}e^A$ does not have an eigenvalue equal to 1. The matrix $Be^{tA}$ is orthogonal for all $t\in\mR$ and thus have eigenvalues on the unit circle in $\mathbb{C}$. The numerics indicate that $\det(I-Be^{tA}) > 0$ for all $t\in[0,1]$ and in all even dimensions $n$. However, the value $t=1$ seems to play a special role: If we define $t_n$ to be the smallest positive number such that $\det(I-Be^{t_n A}) = 0$, we conjecture that the sequence $(t_n)$ is strictly decreasing with $\lim_{n\to\infty} t_n = 1$.

\section{The general case $2\leq p\leq\infty$}

Recall that the Dirichlet problem for $p>2$ is
\begin{align*}
u(x) &= \tfrac{3}{p+1}\cdot\tfrac{1}{2\epsilon}\int_{x-\epsilon}^{x+\epsilon} u(y)\dd y + \tfrac{p-2}{p+1}\cdot\tfrac{1}{2}\Big(u(x-\epsilon) + u(x+\epsilon)\Big), && x\in[-1,1],\\
u(x) &= f(x), && x\in[-1-\epsilon,-1),\\
u(x) &= f(x), && x\in(1,1+\epsilon],
\end{align*}
where $\epsilon = 2/n$.

As usual, we let $x_k := -1 + (k-1)\epsilon$ denote the nodes, and we define the functions $v_k\colon[0,1]\to\mR$ as
\[v_k(t) := u(x_k + \epsilon t),\qquad  k = 0,\dots,n+1.\]
Also, $f_l,f_r\colon[0,1]\to\mR$ are given by
\[f_l(t) := f(x_0 + \epsilon t),\qquad f_r(t) := f(x_{n+1} + \epsilon t).\]

In order to derive the solution, it is instructive to first look at the case $p=\infty$. The equation is then
\[u(x) = \frac{1}{2}\Big(u(x-\epsilon) + u(x+\epsilon)\Big),\]
and for $k=1,\dots,n$ we have
\[v_k(t) = \frac{1}{2}\Big(v_{k-1}(t) + v_{k+1}(t)\Big).\]
Now, $v_0(t) = f_l(t)$ for $0\leq t< 1$ and $v_{n+1}(t) = f_r(t)$ for $0<t\leq 1$. Thus,
\begin{equation}
\bfv(t) = \frac{1}{2}(L+L^\top)\bfv(t) + \frac{1}{2}\big(f_l(t)\bfe_1 + f_r(t)\bfe_n\big)\qquad \text{for $0<t<1$,}
\label{eq:veq}
\end{equation}
where
\[\bfv(t) := [v_1(t),\dots,v_n(t)]^\top\in\mR^n,\qquad L := \begin{bmatrix}
	\bf0^\top & 0\\ I_{n-1} & \bf0
\end{bmatrix}\in\mR^{n\times n}.\]
The well-known matrix $2I-L-L^\top$ is invertible and a direct calculation confirms that $(2I-L-L^\top)\bfw_l = \bfe_1$ and $(2I-L-L^\top)\bfw_r = \bfe_n$ where
\[\bfw_l := \sum_{k=1}^n\left(1 - \frac{k}{n+1}\right)\bfe_k,\qquad \bfw_r := \sum_{k=1}^n\frac{k}{n+1}\bfe_k.\]
The solution of the linear equation \eqref{eq:veq} is therefore
\[\bfv(t) = \sum_{k=1}^n \left[\left(1 - \frac{k}{n+1}\right)f_l(t) + \frac{k}{n+1}f_r(t)\right]\bfe_k,\qquad 0<t<1.\]

\begin{figure}[h]
\centering
\includegraphics{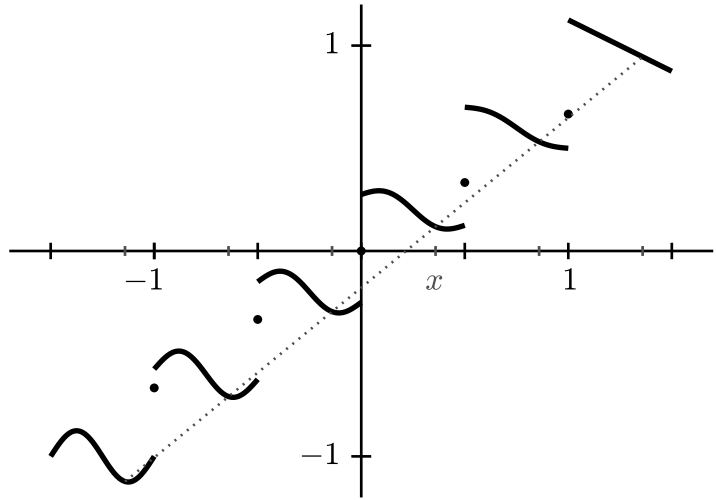}%
\caption{The solution of \eqref{eq:dirp} with $p=\infty$ and $n=4$. The dotted line indicates how $u(x)$ can be constructed from the boundary values.}%
\label{fig_inftynonconst}%
\end{figure}

For $t=0$ we have
\[v_k(0) = \frac{1}{2}
\begin{cases}
f_l(0) + v_{2}(0),\qquad &\text{for $k = 1$,}\\
v_{k-1}(0) + v_{k+1}(0), &\text{for $k = 2,\dots,n$,}\\
v_{n}(0) + f_r(1),\qquad &\text{for $k = n+1$.}
\end{cases}
\]
If we define
\[\hat{\bfv}_0 := [v_1(0),\dots,v_n(0), v_{n+1}(0)]^\top\in\mR^{n+1},\qquad \hat{L} := \begin{bmatrix}
	\bf0^\top & 0\\ I_{n} & \bf0
\end{bmatrix}\in\mR^{(n+1)\times (n+1)},\]
we get the same equation for $\hat{\bfv}_0$ as we did for $\bfv(t)$, but in one dimension higher. Thus,
\[\hat{\bfv}_0 = \sum_{k=1}^{n+1} \left[\left(1 - \frac{k}{n+2}\right)f_l(0) + \frac{k}{n+2}f_r(1)\right]\hat{\bfe}_k.\]

Figure \ref{fig_inftynonconst} shows the graph of the solution to the infinity-equation
\[u(x) = \frac{1}{2}\big(u(x-1/2) + u(x+1/2)\big), \qquad x\in[-1,1],\]
with boundary values
\[u(x) =
\begin{cases}
\frac{1}{8}\sin\big(4\pi(x+3/2)\big) - 1,\qquad &x\in[-3/2,-1),\\
\frac{9}{8} - \frac{x-1}{2}, & x\in(1,3/2].
\end{cases}\]
Notice the jumps $\lim_{x\to x_k^-}u(x) < u(x_k) < \lim_{x\to x_k^+}u(x)$.

We now turn to the case $2\leq p<\infty$.
For $k=0,\dots,n+1$ define the integrals
\[V_k := V_k(1)\qquad\text{where}\qquad V_k(t) := \int_0^t v_k(s)\dd s.\]
Then
\[V_0(t) = \int_0^t f_l(s)\dd s =: F_l(t)\qquad\text{and}\qquad V_{n+1}(t) = \int_0^t f_r(s)\dd s =: F_r(t).\]
The equation can be written as
\begin{align*}
v_k(t) &= u(x_k+\epsilon t)\\
       &= \frac{3}{p+1}\cdot\frac{1}{2\epsilon}\int_{x_{k-1}+\epsilon t}^{x_{k+1}+\epsilon t} u(y)\dd y + \frac{p-2}{p+1}\cdot\frac{1}{2}\Big(u(x_{k-1}+\epsilon t) + u(x_{k+1}+\epsilon t)\Big)\\
         &= \frac{3}{p+1}\cdot\frac{1}{2}\Big( V_{k-1} - V_{k-1}(t) + V_k + V_{k+1}(t)\Big ) + \frac{p-2}{p+1}\cdot\frac{1}{2}\Big( v_{k-1}(t) + v_{k+1}(t)\Big)
\end{align*}
for $k=1,\dots,n$. We have $v_0(t) = f_l(t)$ when $0\leq t< 1$ and $v_{n+1}(t) = f_r(t)$ when $0< t\leq 1$. Thus,
\begin{align*}
\bfV'(t) &= \bfv(t)\\
         &= \frac{3}{p+1}\cdot\frac{1}{2}\left( (L^\top - L)\bfV(t) + F_r(t)\bfe_n - F_l(t)\bfe_1 + (L+I)\bfV + F_l\bfe_1\right)\\
				 &\qquad {} + \frac{p-2}{p+1}\cdot\frac{1}{2}\left( (L+L^\top)\bfV'(t) + f_l(t)\bfe_1 + f_r(t)\bfe_n\right).
\end{align*}
That is,
\begin{equation}
E_p\bfV'(t) = A\bfV(t) + \bfF(t) + \frac{1}{2}(L+I)\bfV,\qquad 0<t<1,
\label{eq:Veq}
\end{equation}
where
\[E_p := \frac{1}{3}\left( (p+1)I - \frac{p-2}{2}(L+L^\top) \right),\]
and where
\[\bfF(t) := \frac{1}{2}\left( \big(F_l(1) - F_l(t)\big)\bfe_1 + F_r(t)\bfe_n \right) + \frac{p-2}{6}\left( f_l(t)\bfe_1+f_r(t)\bfe_n\right).\]
The matrix $E_p$ is invertible since it is diagonal dominant. Since $\bfV(t)$ is continuous with $\bfV(0) = 0$, the solution of \eqref{eq:Veq} is
\begin{equation}
\begin{aligned}
\bfV(t) &= e^{tE_p^{-1}A}\left(0 + \int_0^t e^{-sE_p^{-1}A}E_p^{-1}\left[\bfF(s) + \frac{1}{2}(L+I)\bfV\right]\dd s\right)\\
        &= e^{tE_p^{-1}A}\int_0^t e^{-sE_p^{-1}A}E_p^{-1}\bfF(s)\dd s + \frac{1}{2}\left(e^{tE_p^{-1}A}-I\right)A^{-1}(L+I)\bfV.
\end{aligned}
\label{eq:Vform}
\end{equation}
This gives the linear equation
\begin{equation}
\left[ I - \frac{1}{2}\left(e^{E_p^{-1}A}-I\right)A^{-1}(L+I)\right]\bfV = e^{E_p^{-1}A}\int_0^1 e^{-sE_p^{-1}A}E_p^{-1}\bfF(s)\dd s
\label{eq:lineqV}
\end{equation}
for $\bfV = \bfV(1)$, which, numerically, seems to be nondegenerate.
In fact, one can show that
\[2A\left[ I - \frac{1}{2}\left(e^{E_p^{-1}A}-I\right)A^{-1}(L+I)\right](L^\top + I)^{-1} = I - e^{AE_p^{-1}}\tilde{B}\]
and the question of solvability of \eqref{eq:lineqV} is, at least for $p$ close to 2, equivalent to the solvability for $\bfv_0$ in the previous Section.

The formula for $u$ on the interior of the intervals $(x_k,x_{k+1})$ follows now from \eqref{eq:Veq}:
\[\bfv(t) = E_p^{-1}\left(A\bfV(t) + \bfF(t) + \frac{1}{2}(L+I)\bfV\right),\qquad 0<t<1.\]
When $t=0$ we have
\[v_k(0) = \frac{3}{p+1}\cdot\frac{1}{2}\Big( V_{k-1}  + V_k\Big ) + \frac{p-2}{p+1}\cdot\frac{1}{2}\Big( v_{k-1}(0) + v_{k+1}(0)\Big)\]
for $k=1,\dots,n$. Now,
\[v_1(0) = \frac{3}{p+1}\cdot\frac{1}{2}\Big( F_l(1) + V_1\Big ) + \frac{p-2}{p+1}\cdot\frac{1}{2}\Big( f_l(0) + v_{2}(0)\Big),\]
and
\[v_{n+1}(0) = \frac{3}{p+1}\cdot\frac{1}{2}\Big( V_{n} + F_r(1)\Big ) + \frac{p-2}{p+1}\cdot\frac{1}{2}\Big( v_{n}(0) + f_r(1)\Big).\]
We set 
\[Q := \begin{bmatrix}
	I\\ \bf0^\top
\end{bmatrix}\in \mR^{(n+1)\times n},\]
and get the solvable equation
\[\hat{E}_{p}\hat{\bfv}_0 = \frac{1}{2}\Big( \hat{L} + \hat{I}\Big )Q\bfV + \frac{1}{2}\big(F_l(1)\hat{\bfe}_1 + F_r(1)\hat{\bfe}_{n+1}\big) + \frac{p-2}{6}\big(f_l(0)\hat{\bfe}_1  + f_r(1)\hat{\bfe}_{n+1}\big)\]
for $\hat{\bfv}_0 := [v_1(0),\dots, v_{n+1}(0)]^\top$ in $\mR^{n+1}$.

\begin{figure}[h]
\centering
\includegraphics{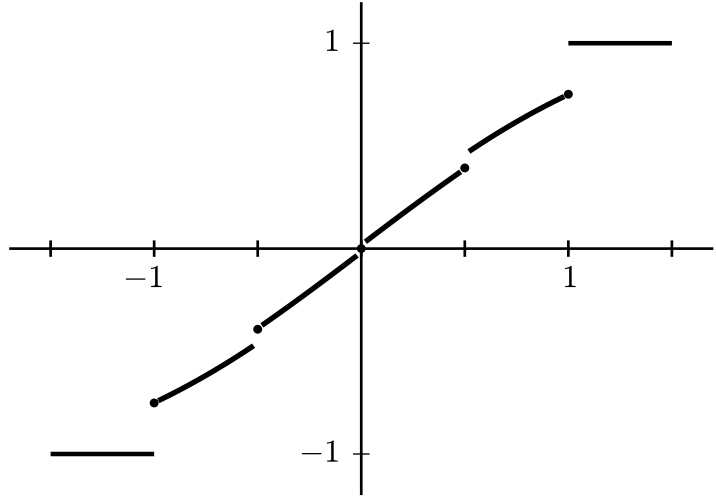}%
\caption{The solution of \eqref{eq:dirp} with $p=5$, $n=4$, and boundary data $\pm 1$.}%
\label{figp4}%
\end{figure}

Some simplifications can be made when the boundary values are constant. If $f_l(t) = -1$ and $f_r(t) = 1$, then
\begin{align*}
\bfF(t) &= \frac{1}{2}\left( \big(F_l- F_l(t)\big)\bfe_1 + F_r(t)\bfe_n \right) + \frac{p-2}{6}\left( f_l(t)\bfe_1+f_r(t)\bfe_n\right)\\
        &= \frac{t}{2}( \bfe_1 + \bfe_n) + \frac{p-2}{6}\bfe_n - \frac{p+1}{6}\bfe_1\\
				&= \frac{t}{2}( \bfe_1 + \bfe_n) + \frac{p-2}{6}(\bfe_n-\bfe_1) + \frac{1}{2}(L+I)\tilde{\bfe},
\end{align*}
and an integration by parts,
\begin{align*}
\int_0^t se^{-sE_p^{-1}A}E_p^{-1}\dd s &= -\bigg|_0^t se^{-sE_p^{-1}A}(E_p^{-1}A)^{-1}E_p^{-1} + \int_0^t e^{-sE_p^{-1}A}(E_p^{-1}A)^{-1}E_p^{-1}\dd s\\
	&= - te^{-tE_p^{-1}A}A^{-1} - \left(e^{-tE_p^{-1}A} - I\right)A^{-1}E_p A^{-1},
\end{align*}
yields
\begin{align*}
\bfV(t) &= e^{tE_p^{-1}A}\int_0^t e^{-sE_p^{-1}A}E_p^{-1}\Big[\tfrac{s}{2}( \bfe_1 + \bfe_n) + \tfrac{p-2}{6}( \bfe_n - \bfe_1) + \tfrac{1}{2}(L+I)(\bfV+\tilde{\bfe})\Big]\dd s\\
        &= e^{tE_p^{-1}A}\left(- te^{-tE_p^{-1}A}A^{-1} - \left(e^{-tE_p^{-1}A} - I\right)A^{-1}E_p A^{-1}\right)\frac{1}{2}( \bfe_1 + \bfe_n)\\
				&\qquad {} + e^{tE_p^{-1}A}\int_0^t e^{-sE_p^{-1}A}E_p^{-1}\left[\frac{p-2}{6}(\bfe_n - \bfe_1) + \frac{1}{2}(L+I)\big(\bfV+\tilde{\bfe}\big)\right]\dd s\\
				&= -\left( tI + \left(I - e^{tE_p^{-1}A}\right)A^{-1}E_p \right)\tilde{\bfe}\\
				&\qquad {} - \left(I - e^{tE_p^{-1}A}\right)A^{-1}\left( \frac{p-2}{6}( \bfe_n - \bfe_1) + \frac{1}{2}(L+I)\big(\bfV+\tilde{\bfe}\big)\right)\\
				&= \left(e^{tE_p^{-1}A} - I\right)A^{-1}\left( \frac{p-2}{6}( \bfe_n - \bfe_1) + \frac{1}{2}(L+I)\big(\bfV+\tilde{\bfe}\big) + E_p\tilde{\bfe}\right) - t\tilde{\bfe}.
\end{align*}
The equation for $\bfV + \tilde{\bfe}$ is then
\begin{align*}
\left[ I - \tfrac{1}{2}\left(e^{E_p^{-1}A}-I\right)A^{-1}(L+I)\right]\big(\bfV + \tilde{\bfe}\big)
	&= \left(e^{E_p^{-1}A} - I\right)A^{-1}\left( \tfrac{p-2}{6}( \bfe_n - \bfe_1) + E_p\tilde{\bfe} \right),
\end{align*}
and a differentiation gives
\[\bfv(t) = \bfV'(t) = e^{tE_p^{-1}A}\left( \frac{p-2}{6}E_p^{-1}( \bfe_n - \bfe_1) + \frac{1}{2}E_p^{-1}(L+I)\big(\bfV + \tilde{\bfe}\big) + \tilde{\bfe} \right) - \tilde{\bfe}.\]

\begin{figure}[h]
\centering
\includegraphics{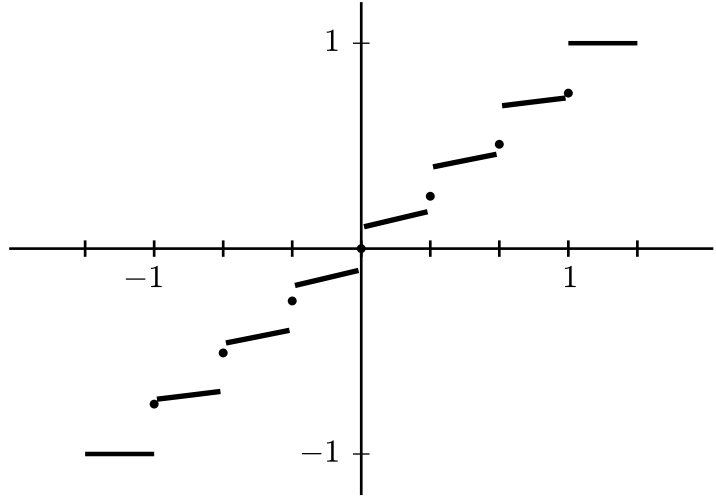}%
\caption{$p=100$ and $n=6$. Constant boundary $\pm 1$.}%
\label{figp100}%
\end{figure}

Finally, the equation for the values $\hat{\bfv}_0 := [v_1(0),\dots, v_{n+1}(0)]^\top$ at the nodes becomes
\begin{align*}
\hat{E}_{p}\hat{\bfv}_0
	&= \frac{1}{2}\Big( \hat{L} + \hat{I}\Big )Q\bfV + \frac{1}{2}\big(F_l(1)\hat{\bfe}_1 + F_r(1)\hat{\bfe}_{n+1}\big) + \frac{p-2}{6}\big(f_l(0)\hat{\bfe}_1  + f_r(1)\hat{\bfe}_{n+1}\big)\\
	&= \frac{1}{2}\Big( \hat{L} + \hat{I}\Big )Q\bfV +  \frac{p+1}{6}\big(\hat{\bfe}_{n+1} - \hat{\bfe}_1 \big).
\end{align*}

\begin{figure}[h]
\centering
\includegraphics{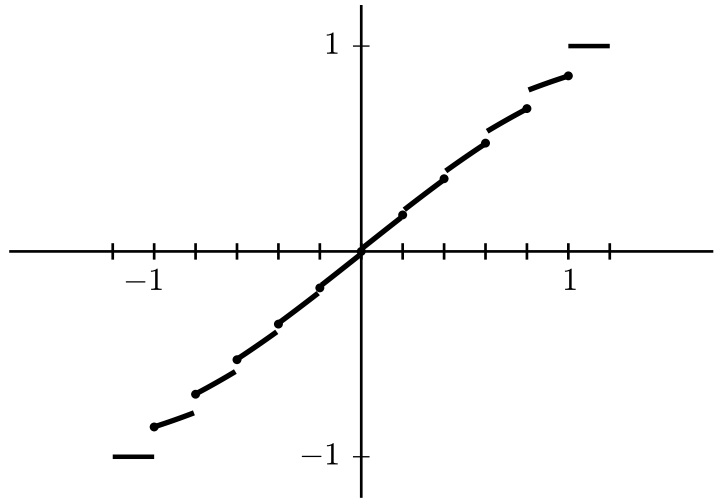}%
\caption{$p=25$ and $n=10$. Constant boundary $\pm 1$.}%
\label{fign10p25}%
\end{figure}

The Figures \ref{figp4} - \ref{fign10p25} show the graph of the solution to \eqref{eq:dirp} for various $n$ and $p$, and when the boundary values are constant. As $p\to\infty$ the solution converges (slowly) to the solution of the infinity-equation.

When $p>2$, the matrix $E_p^{-1}A$ is not skew-symmetric and $e^{tE_p^{-1}A}$ is no longer orthogonal. However, since $E_p$ is symmetric and positive definite, the eigenvalues of $E_p^{-1}A$ are still purely imaginary and the eigenvalues of $e^{tE_p^{-1}A}$ are again on the form $\cos(\lambda t) + i\sin(\lambda t)$. The mean value solutions are therefore piecewise trigonometric also in the case $2<p<\infty$.

\bibliographystyle{alpha}
\bibliography{C:/Users/Karl_K/Documents/PhD/references}

\end{document}